\documentclass[12pt]{amsart}

\usepackage{amsmath,amssymb,amsthm,graphicx,mathrsfs,url}
\newcommand{\C}{\mathbb{C}}

\newtheorem{theorem}{Theorem}

\title[Resonance counting function for surface with cusps]{A note on the resonance counting function for surfaces with cusps}
\author{Yannick Bonthonneau}
\email{yannick.bonthonneau@ens.fr}
\address{DMA, U.M.R. 8553 CNRS, \'Ecole Normale Superieure, 45 rue d'Ulm,
75230 Paris cedex 05, France}
\begin{document}

\begin{abstract}
We prove sharp upper bounds for the number of resonances in boxes of size $1$ at high frequency 
for the Laplacian on finite volume surfaces with hyperbolic cusps. 
As a corollary, we obtain a Weyl asymptotic for the number of resonances in balls of size $T\to \infty$ with remainder $O(T^{3/2})$.
\end{abstract}

\maketitle

In this short note, we intend to prove sharp bounds on resonance-counting functions for the Laplacian on finite volume surfaces with hyperbolic cusps. Let $M$ be a complete non-compact surface, equipped with a Riemannian metric $g$. We assume that $(M,g)$ can be decomposed as the union of a compact manifold with boundary and a finite number of hyperbolic cusps, each one being isometric to
\begin{equation*}
(a,+\infty)_y\times \mathbb{S}_\theta^1 \textrm{ with metric }\frac{dy^2+d\theta^2}{y^2}
\end{equation*}
for some $a>0$. The spectral properties of the Laplacian $\Delta_g$ were first studied by  Selberg \cite{MR1117906} and Lax-Phillips \cite{MR0562288} in constant negative curvature, and by Colin-de-Verdi\`ere \cite{MR699488}, M\"uller \cite{MR1172692}, Parnovski \cite{MR1348800} in the  non-constant curvature setting. \\

On such surfaces, the resolvent $R(s)=(\Delta_g-s(1-s))^{-1}$ of the Laplacian admits a meromorphic extension from $\{\Re s>1/2\}$ to $\mathbb{C}$ as an operator mapping $L^2_{\rm comp}$ to $L^2_{\rm loc}$ and the natural discrete spectral set  for $\Delta_g$ is the set of poles denoted by
\begin{equation*}
\mathcal{R} \subset \{ s \in \mathbb{C} | \,\, \Re s \leq 1/2 \}\cup (1/2, 1].
\end{equation*}
The poles are called \emph{resonances} and are counted with multiplicity $m(s)$ (the multiplicity $m(s)$ is defined below and corresponds, for all but finitely many resonances, to the rank of the residue of the resolvent at $s$). We shall recall in the next section how the set of resonances is built. To study their distribution in the complex plane, we define two counting functions :
\begin{align}
N_{\mathcal{R}} (T) &:= \sum_{s\in \mathcal{R}, |s-1/2| \leq T} m(s) \\
N_{\mathcal{R}}(T,\delta) &:= \sum_{s\in \mathcal{R}, |s-1/2 - \imath T| \leq \delta} m(s).
\end{align}

The first result on the resonance counting function was proved by Selberg \cite[p. 25]{MR1117906} for the special case of 
hyperbolic surfaces with finite volume : the following Weyl type asymptotic expansion holds as $T\to \infty$
\begin{equation}\label{Selberg}
N_{\mathcal{R}}(T) = \frac{{\rm Vol}(M)}{2\pi} T^2+ C_0T\log(T)+C_1T+O\Big(\frac{T}{\log(T)}\Big)
\end{equation}
for some explicit constants $C_0,C_1$.  
In variable curvature, M\"uller gives a Weyl asymptotic  \cite[Th. 1.3.a]{MR1172692} of the form 
\begin{equation*}
N_{\mathcal{R}}(T)=\frac{{\rm Vol}(M)}{2\pi}T^2+o(T^2)
\end{equation*}
and this was improved by Parnovski \cite{MR1348800} who showed that for all $\epsilon>0$
\begin{equation}
N_{\mathcal{R}}(T) = \frac{{\rm Vol}(M)}{2\pi} T^2 + O(T^{3/2+\epsilon}).
\end{equation}
Parnovski's proof relies on a Weyl type asymptotic expansion involving the scattering phase $\mathcal{S}(T)$ (see next section for a precise definition) :
\begin{equation}\label{eq:Phase_estimate}
2\pi N_d(T) + \mathcal{S}(T) = \frac{{\rm Vol}(M)}{2}T^2 - 2k T \ln T + O(T),
\end{equation}
where $k$ is the number of cusps, and $N_d$ is the counting function for the $L^2$ eigenvalues of $\Delta_g$ embedded in the continuous spectrum.

Using a Poisson formula proved by M\"uller \cite{MR1172692} and estimate \eqref{eq:Phase_estimate}, we are able to improve the results of Parnovski :
\begin{theorem}
For $T>1$, and $0\leq \delta \leq T/2$, the following estimates hold
\begin{align}
N_{\mathcal{R}}(T,\delta) &= O(T \delta + T) \label{eq:small_box}, \\
N_{\mathcal{R}} (T) &= \frac{{\rm Vol}(M)}{2\pi} T^2 + O(T^{3/2}).
\end{align}
\end{theorem}

In the first estimate with $\delta=1$, the exponent in $T$ is sharp in general, as can be seen from Selberg's result \eqref{Selberg}
which implies that there is $C>0$ such that as $T\to \infty$
\begin{equation}
N(T,1) = C T + O\left(\frac{T}{\log T}\right)
\end{equation}

In $n$-dimensional Euclidan scattering, upper bounds $O(T^{n-1})$ on the number of resonances in boxes of fixed size 
at frequency $T$ were obtained by Petkov-Zworski 
\cite{MR1704278} using Breit-Wigner approximation and the scattering phase ; our scheme of proof is inspired from their approach. 
Their result was extended to the case of non-compact perturbations of the Laplacian by Bony \cite{MR1853138}. 
In general, it is expected that the number of resonances in such boxes is controlled by the (fractal) dimension of the trapped set (see for example Zworski \cite{MR1688441}, Guillop\'e-Lin-Zworski \cite{MR2036371}, Sj\"ostrand-Zworski \cite{MR2309150}, Datchev-Dyatlov \cite{MR3077910}).\\

\textbf{Acknowledgement}. We thank M. Zworski for his suggestion which shortened significantly the argument of proof. We also thank J-F. Bony for sending us his work, and Colin Guillarmou and Nalini Anantharaman for their fruitful advice.

\section{Preliminaries}

We start by recalling well-known facts on scattering theory on surfaces with cusps, and we refer to the article of M\"uller \cite{MR1172692} for details. Let  $(M,g)$ be a complete Riemannian surface that can be decomposed as follows: 
\begin{equation*}
M = M_0 \cup Z_1 \cup \dots \cup Z_k,
\end{equation*}
where $M_0$ is a compact surface with smooth boundary, and $Z_j$ are hyperbolic cusps
\begin{equation*}
Z_j \simeq (a_j, +\infty ) \times \mathbb{S}^1, \quad j= 1 \dots k,
\end{equation*}
with $a_j > 0$ and the metric on $Z_j$ in coordinates $( y,\theta) \in  (a_j , + \infty) \times \mathbb{S}^1$ is
\begin{equation*}
{d}s^2 = \frac{\mathrm{d}y^2 + \mathrm{d}\theta^2}{y^2}.
\end{equation*}
Notice that the surface has finite volume when equipped with this metric.

The non-negative Laplacian $\Delta$ acting on $C_0^\infty(M)$ functions has a unique self-adjoint extension to $L^2(M)$ and its spectrum consists of
\begin{enumerate}
\item Absolutely continuous spectrum $\sigma_{ac}=[1/4, +\infty)$ with multiplicity $k$ (the number of cusps).
\item Discrete spectrum $\sigma_d =\{\lambda_0 = 0 < \lambda_1 \leq \dots \leq \lambda_i \leq \dots \}$, possibly finite, and which may contain embedded eigenvalues in the continuous spectrum. To $\lambda \in \sigma_d$, we associate a family of orthogonal eigenfunctions that generate its eigenspace $(u_\lambda^{i})_{i=1 \dots d_\lambda} \in L^2(M) \cap C^\infty(M)$.
\end{enumerate}

The generalized eigenfunctions associated to the absolutely continuous spectrum are the Eisenstein functions, $(E_j(x,s))_{i=1 \dots k}$. Each $E_j$ is a meromorphic family (in $s$) of smooth functions on $M$. Its poles are contained in the open half-plane $\{\Re s < 1/2\}$ or in $(1/2, 1]$. The Eisenstein functions are characterized by two properties : 
\begin{enumerate}
\item $\Delta_g E_j(.,s) = s(1-s)E_j(.,s)$ 
\item In the cusp $Z_i$, $i=1 \dots k$, the zeroth Fourier coefficient of $E_j$ in the $\theta$ variable 
equals $\delta_{ij} y_i^s + \phi_{ij}(s) y_i^{1-s}$ where $y_i$ denotes the $y$ coordinate in the cusp $Z_i$ and $\phi_{ij}(s)$ is a meromorphic function of $s$.
\end{enumerate}

We can collect the \emph{scattering coefficients} $\phi_{ij}$ in a meromorphic family of matrices, $\phi(s) = (\phi_{ij})_{ij}$ called \emph{scattering matrix}. We denote its determinant by $\varphi(s) = \det \phi(s)$. Then the following identities hold
\begin{equation*}
\phi(s)\phi(1-s)=Id, \quad \overline{\phi(s)} = \phi(\overline{s}), \quad \phi(s)^\ast = \phi(\overline{s}).
\end{equation*}
The line $\Re s =1/2$ corresponds to the continuous spectrum. On that line, $\phi(s)$ is unitary, $\varphi(s)$ has modulus $1$. We also define the scattering phase 
\begin{equation}
\mathcal{S}(T) = - \int_0^T \frac{\varphi'}{\varphi}(\frac{1}{2} + \imath t) \mathrm{d} t
\end{equation}

The set of poles of $\varphi$, $\phi$ and $(E_j)_{j=1\dots k}$ is the same, we call them them \emph{scattering poles} and we shall denote $\Lambda$ this set. It is contained in $\{\Re s < 1/2\} \cup (1/2, 1]$. The union of this set with the set of $s\in \C$ such that $s(1-s)$ is an $L^2$ eigenvalue, is called the resonance set, and denoted $\mathcal{R}$. Following \cite[pp.287]{MR1172692}, the multiplicities $m(s)$ are defined as :
\begin{enumerate}
	\item If $\Re s \geq 1/2$, $s \neq 1/2$, $m(s)$ is the dimension of $\ker_{L^2}(\Delta_g-s(1-s))$.
	\item If $\Re s < 1/2$, $m(s)$ is the dimension of $\ker_{L^2}(\Delta_g-s(1-s))$ minus the order of $\varphi$ at $s$.
	\item $m(1/2)$ equals $({\rm Tr}(\phi(1/2)) + k)/2$ plus twice the dimension of $\ker_{L^2}(\Delta_g-1/4)$.
\end{enumerate}

For convenience, we define two counting functions for the discrete spectrum and the poles of $\varphi$: 
\begin{align}
N_d(T) & := \sum_{|s_i-1/2| \leq T}m(s_i), \\
N_{\Lambda}(T) & : = \sum_{s\in \Lambda, |s-1/2| \leq T}m(s),
\end{align}
so that
\[N_{\mathcal{R}}(T):=  \sum_{s\in \mathcal{R}, |s-1/2| \leq T}m(s)= 2 N_d(T) + N_{\Lambda}(T).\]

\section{Main observation}

In this Section, we explain how to obtain estimate for $N_{\mathcal{R}}(T)$ in boxes at high frequency. 

From the asymptotic expansion \eqref{eq:Phase_estimate}, we deduce that for $0\leq \delta \leq T/2$,
\begin{equation}\label{parnov}
2\pi(N_d(T+\delta) - N_d(T-\delta)) + \mathcal{S}(T+ \delta) -\mathcal{S}(T-\delta) = 2{\rm Vol}(M)T\delta  - 4 k \delta \ln T + O(T).
\end{equation}
Next, we recall the Poisson formula for resonances proved by M\"uller \cite[Th. 3.32]{MR1172692}
\begin{equation}\label{poisson}
\mathcal{S}'(T) = \log \frac{1}{q} + \sum_{\rho \in \Lambda} \frac{1-2\Re \rho}{(\Re \rho -1/2)^2 + (\Im \rho - T)^2}.
\end{equation}
where $q$ is some constant (not necessarily $<1$). Let $C>1$, $0<\epsilon<1$ and 
\begin{equation*}
\Omega_{T,\delta} := \{ s \in \mathbb{C}; \,\, |s-1/2 - \imath T| \leq \delta / C \text{ and } 0 \leq 1/2 - \Re s \leq \epsilon \delta \}.
\end{equation*}

Then, for $s\in \Omega_{T,\delta}$, 
\begin{equation*}
\int_{[T-\delta,T+\delta]} \frac{1-2\Re s}{(\Re s -1/2)^2 + (t - \Im s)^2}\mathrm{d}t = 2\left[ \arctan \frac{t - \Im s}{1/2 - \Re s} \right]_{T-\delta}^{T+\delta}
\end{equation*}
The addition formula for $\arctan$, with $x, y> 0$ and $xy > 1$ is given by 
\begin{equation*}
\arctan x + \arctan y = \pi + \arctan \frac{x+y}{1-xy} 
\end{equation*}
thus
\begin{equation*}
\begin{split}
\int_{[T-\delta,T+\delta]} \frac{1-2\Re s}{(\Re s -1/2)^2 + (t - \Im s)^2}\mathrm{d}t &=  2\pi - 2 \arctan \frac{2\delta (1/2-\Re s)}{\delta^2 - |s-1/2 -\imath T|^2}\\
&\geq 2\pi - 2 \arctan \tilde{C}\epsilon,
\end{split}
\end{equation*}
where $\tilde{C}$ is set to be $2/(1-1/C^2)$. For $\epsilon$ small enough, this is bigger than, say $\pi$.

Since all but a finite number of terms in \eqref{poisson} are positive, we have :
\begin{equation*}
\mathcal{S}(T+ \delta) -\mathcal{S}(T-\delta) \geq O(\delta) + \sum_{\rho \in \Lambda \cap \Omega_{T,\delta}} \pi.
\end{equation*}

Combining with \eqref{parnov}, we deduce that 
\begin{equation*}
N_d(T+\delta)-N_d(T-\delta) + \# \Lambda \cap \Omega_{T,\delta} = O( T \delta) + O(T) + O(\delta).
\end{equation*}

This is the content of \eqref{eq:small_box} in our main theorem.

\section{Consequence}

Now, we proceed to prove the second part of our theorem. We will follow the method of M\"uller \cite[pp. 282]{MR1172692}, 
which is a global and quantitative version of the argument used in the previous section. Integrating the Poisson formula over $[-T,T]$, we relate the scattering phase asymptotics to the poles of $\phi$. Using the $\arctan$ addition formula, we are left with the sum of $N_\Lambda(T)$ and an expression with $\arctan$'s (equation (4.9) in \cite{MR1172692}) :
\begin{equation}
\begin{split}
&\frac{1}{2\pi} \mathcal{S}(T) = \\
\frac{1}{2} N_\Lambda(T)& + \frac{1}{2\pi} \sum_{\rho \in \Lambda, \Re \rho < 1/2} \arctan\left[ \frac{1-2\Re \rho}{|\rho - 1/2|^2} T \left(1 - \frac{T^2}{|\rho - 1/2 |^2} \right)^{-1} \right] + O(T).
\end{split}
\end{equation}

The sum is then split between $\{1\}$ the poles in $\{ |T- |\rho -1/2|| > T^{1/2} \}$, and $\{2\}$, the others. M\"uller proved that the sum $\{1\}$ is $O(T^{3/2})$. The sum $\{2\}$ can be bounded by 
\begin{equation*}
\frac{1}{4}( N_\Lambda(T+ \sqrt{T}) - N_\Lambda(T- \sqrt{T})).
\end{equation*}

From \cite[Cor. 3.29]{MR1172692}, we also recall that 
\begin{equation*}
\sum_{\eta \in \Lambda, \eta \neq 1/2} m(\eta)\frac{1 - 2\Re \eta}{|\eta - 1/2|^2} < \infty.
\end{equation*}
Consider the set $\tilde{\Lambda}= \{\eta \in \Lambda; \: (2\Re \eta - 1)^2 >  \Im \eta, \: |\eta| > 1 \}$. 
On $\tilde{\Lambda}$, we have that $|\eta-1/2|^{1/2} \leq 1-2\Re \eta$, thus  
\begin{equation*}
\sum_{\eta \in \tilde{\Lambda}, \eta \neq 1/2} m(\eta) \frac{1}{|\eta-1/2|^{3/2}} < \infty.
\end{equation*}
If $\tilde{n}(T)$ is the counting function for $\tilde{\Lambda}$, we deduce that
\begin{equation*}
\sum_{k=1}^\infty \tilde{n}(k) \left[\frac{1}{k^{3/2}} - \frac{1}{(k+1)^{3/2}}\right] < \infty.
\end{equation*}
Since $\tilde{n}$ is non-decreasing, $\tilde{n}(k) = o(k^{3/2})$. Now, 
\begin{equation*}
N_\Lambda(T-\sqrt{T}) - N_\Lambda(T+ \sqrt{T}) \leq \tilde{n}(T) + N_{\mathcal{R}}( T,\sqrt{T}) + N(T,\sqrt{T}).
\end{equation*}
This concludes the proof.

\bibliographystyle{alpha}
\bibliography{../../bibliographies/biblio_article_Counting_Resonances}

\end{document}